\documentclass[a4paper,reqno]{amsart}

\usepackage{amsthm}
\usepackage{amssymb}
\usepackage{amsmath}
\usepackage{mathtools}
\usepackage{pdfpages}
\usepackage{tikz-cd}
\usepackage{exscale,cite,color,amsopn,fullpage}
\usepackage[colorlinks=true,pdftex,unicode=true,linktocpage,bookmarksopen,hypertexnames=false]{hyperref}

\newcommand{\free}[1]{\langle#1\rangle}

\renewcommand{\leq}{\leqslant}
\renewcommand{\geq}{\geqslant}

\DeclareMathOperator{\Cen}{Z}
\DeclareMathOperator{\clK}{clKdim}
\DeclareMathOperator{\Ext}{Ext}

\DeclareMathOperator{\GK}{GKdim}

\DeclareMathOperator{\id}{id}
\DeclareMathOperator{\Map}{Map}
\DeclareMathOperator{\Ret}{Ret}
\DeclareMathOperator{\rk}{rk}

\DeclareMathOperator{\Span}{Span}
\DeclareMathOperator{\Sym}{Sym}

\makeatletter

\numberwithin{equation}{section}

\newtheorem{thm}{Theorem}[section]
\newtheorem{lem}[thm]{Lemma}
\newtheorem{cor}[thm]{Corollary}
\newtheorem{pro}[thm]{Proposition}
\newtheorem*{mthm}{Theorem}

\theoremstyle{definition}
\newtheorem{dfn}[thm]{Definition}
\newtheorem{exa}[thm]{Example}

\definecolor{Green}{rgb}{0.13, 0.55, 0.13}

\makeatother

\title{Set-theoretic solutions of the Pentagon Equation}
\author{Ilaria Colazzo \and Eric Jespers \and {\L}ukasz Kubat}

\address[I. Collazzo, E. Jespers, {\L}. Kubat]{Department of Mathematics, Vrije Universiteit Brussel, Pleinlaan 2, 1050 Brussel}
\email{ilaria.colazzo@vub.be}
\email{eric.jespers@vub.be}
\email{lukasz.kubat@vub.be}

\thanks{The third author is supported by Fonds voor Wetenschappelijk Onderzoek (Flanders), grant G016117.
	The first and second authors are supported in part by Onderzoeksraad of Vrije Universiteit Brussel and Fonds voor Wetenschappelijk Onderzoek (Flanders), grant G016117.}

\subjclass[2010]{Primary: 16T25, 16T20; Secondary: 16S36, 16S37, 20M25}
\keywords{pentagon equation, set-theoretic solution, structure algebra, Yang--Baxter equation}

\begin{document}

\maketitle

\begin{abstract}
    A set-theoretic solution of the Pentagon Equation on a non-empty set $S$ is a map $s\colon S^2\to S^2$
    such that $s_{23}s_{13}s_{12}=s_{12}s_{23}$, where $s_{12}=s\times{\id}$, $s_{23}={\id}\times s$ and
    $s_{13}=(\tau\times{\id})({\id}\times s)(\tau\times{\id})$ are mappings from $S^3$ to itself and
    $\tau\colon S^2\to S^2$ is the flip map, i.e., $\tau (x,y) =(y,x)$. We give a description of all
    involutive solutions, i.e., $s^2=\id$. It is shown that such solutions are determined by a factorization
    of $S$ as direct product $X\times A \times G$ and a map $\sigma\colon A\to\Sym(X)$, where $X$ is a non-empty
    set and $A,G$ are elementary abelian $2$-groups. Isomorphic solutions are determined by the cardinalities of
    $A$, $G$ and $X$, i.e., the map $\sigma$ is irrelevant. In particular, if $S$ is finite of cardinality
    $2^n(2m+1)$ for some $n,m\geq 0$ then, on $S$, there are precisely $\binom{n+2}{2}$ non-isomorphic
    solutions of the Pentagon Equation.
\end{abstract}

\section*{Introduction}

The Quantum Yang--Baxter Equation plays a crucial role in the theory of two-dimensional integrable systems of field theory
and statistical mechanics. This appeared first in the works of Yang \cite{Yang1967} and Baxter \cite{Baxter1971}.
It also lead to the theory of quantum groups and related areas \cite{Kassel1995}. Recall that a solution of the Quantum
Yang--Baxter Equation is a linear map $R\colon V\otimes V\to V\otimes  V$, where $V$ is a vector space, such that
$R_{12}R_{13}R_{23}= R_{23}R_{13}R_{12}$, where $R_{ij}$ denotes the map $V\otimes V\otimes V \to V\otimes  V\otimes V$
acting as $R$ on the $(i,j)$-th tensor factors and as the identity on the remaining factor. 
Zamolodchikov \cite{Zamol1980,Zamol1981} described a generalization of this equation, the Tetrahedron Equation,
for three-dimensional integrable systems. Maillet in \cite{Ma94} shows that solutions of the Pentagon Equation lead
to solutions of the Tetrahedron Equation.  Recall that a linear map $S\colon V\otimes V\to V\otimes V$ is a solution of
the Pentagon Equation if, on $V\otimes V \otimes V$, one has $S_{23}S_{13}S_{12}=S_{12}S_{23}$. (In the literature
there is no agreement because some authors call $S_{12}S_{13}S_{23}=S_{23}S_{12}$ the Pentagon Equation. However, solutions
of both equations are in bijection via $S\mapsto\tau S\tau$, where $\tau\colon V\otimes V\to V\otimes V$ is the flip map
defined as $\tau(u\otimes v)=v\otimes u$.) Solutions of the Pentagon Equation appear now in various contexts. For example,
Kashaev in \cite{Ka96} has shown that the Pentagon Equation plays the same role for the Heisenberg double as the Quantum
Yang--Baxter Equation does for the Drinfeld double. Militaru \cite{Mi04} has shown that any finite-dimensional Hopf algebra
is characterized by an invertible solution of the Pentagon Equation. For other contexts we refer the reader to \cite{DiMu15}.
Solutions of the Pentagon Equation also appear with different terminologies. For example, for a Hilbert space $H$, a unitary
operator on $H\otimes H$ is said to be multiplicative if it is a solution of the Pentagon Equation \cite{BaSk93}; and in \cite{St98},
for a fixed braided monoidal category $\mathcal{V}$, an arrow $V\colon A\otimes A \to A \otimes A$ in $\mathcal{V}$
is said to be a fusion operator if it satisfies the Pentagon Equation.

Note that if $S$ is a basis of the vector space $V$ then a map $s\colon S\times S\to S\times S$ such that
$s_{23}s_{13}s_{12}=s_{12}s_{23}$ (respectively, $s_{12}s_{13}s_{23}=s_{23}s_{13}s_{12}$), where $s_{ij}$ denotes the map
$S\times S\times S\to S\times S\times S$ acting as $s$ on the $(i,j)$-th components and as the identity on the remaining component,
induces a solution of the Pentagon Equation (respectively, the Quantum Yang--Baxter Equation). In this case, one says that $(S,s)$
is a set-theoretic solution of the Pentagon Equation (respectively, the Quantum Yang--Baxter Equation). Drinfeld,
in \cite{Drinfeld1992}, posed the question of finding these set-theoretic solutions.

In the case of the Yang--Baxter Equation, a subclass of this type of solutions, the non-degenerate involutive solutions, has received
a lot of attention in past two decades; we give only a few references \cite{ESS,CJO2010,CJO2014,GIVdB,Lu2000,Rump2005}. This class of
solutions is not only studied for the applications of the Yang--Baxter Equation in physics, but also for its connection with other topics
in mathematics of recent interest, such as semigroups of $I$-type and Bieberbach groups, bijective $1$-cocycles, radical rings, triply
factorized groups, constructions of semisimple minimal triangular Hopf algebras, regular subgroups of the holomorph and Hopf--Galois extensions,
and groups of central type. Also general non-degenerate bijective solutions have received a lot of attention, see for example \cite{GuVe2017}.
It has been shown \cite{BCJ2016,Bach2018} that a description of all finite bijective non-degenerate solutions of the Yang--Baxter Equation
follows from a description of the algebraic structures called braces and skew braces, as introduced by Rump, and Guarnieri and Vendramin,
respectively. However the classification of these structures, and hence the construction of all finite bijective (involutive) non-degenerate 
set-theoretic solutions  of the Yang--Baxter Equation is presently beyond reach.

Also the set-theoretic solutions of the Pentagon Equation received a large interest; we refer to, for example,
\cite{KaRe07,BaSk93,BaSk03,Jiang2005,Ka2011,KaSe98}. Recently, in \cite{CMM19x}, Catino, Mazzotta and Miccoli
dealt with a restrictive case of set-theoretic solutions of the Pentagon Equation. Namely, they described
all, not necessarily bijective, solutions $(S,s)$ of the Pentagon Equation of the form $s(x,y)=(x\cdot y,x*y)$,
where either $(S,\cdotp)$ or $(S,*)$ is a group. Furthermore, Catino, Mazotta and Stefanelli in \cite{CMS2019}
provided a method to obtain solutions of the Pentagon Equation on the matched product of two semigroups, that is a
semigroup including the classical Zappa--Sz\'ep product.

In this paper we continue the investigations on the description of bijective set-theoretic solutions of the Pentagon Equation. Surprisingly,
we are able to describe all such involutive solutions. This  is in strong contrast with the knowledge for involutive solutions
of the Yang--Baxter Equation.

\begin{mthm}
	Assume that $(S,s)$ is an involutive solution of the Pentagon Equation. Then there exist elementary abelian $2$-groups
	$(A,+)$ and $(G,\cdotp)$, a non-empty set $X$ and a map $\sigma\colon A\to \Sym(X)$ (not necessarily a group morphism)
	such that the set $S$ may be identified with $X\times A\times G$ and then
	\[s((x,a,g),(y,b,h))=((x,a,g\cdot h),(\sigma_{a+b}\sigma_b^{-1}(y),a+b,h))\] for all $x,y\in X$, $a,b\in A$ and $g,h\in G$.
	
	Moreover, all involutive solutions, up to isomorphism, of the Pentagon Equation defined on a set $S$ are in a bijective correspondence
	with decompositions of $S$ as a product $X\times A\times G$, where $X$ is a non-empty set and $A,G$ are elementary abelian $2$-groups.
	In particular, if $S$ is finite of cardinality $|S|=2^n(2m+1)$ for some $n,m\geq 0$ then there exist, up to isomorphism, exactly $\binom{n+2}{2}$
	involutive solutions of the Pentagon Equation defined on $S$.
\end{mthm}

\section{Definitions and examples}

Let $S$ be a non-empty set. A map $s\colon S\times S\to S \times S$ (or a pair $(S,s)$) is said to be a
\emph{set-theoretic solution of the Pentagon Equation} provided \[s_{23}s_{13}s_{12}=s_{12}s_{23},\]
where $s_{12}=s\times{\id}$, $s_{23}={\id}\times s$ and $s_{13}=(\tau\times{\id})({\id}\times s)(\tau\times{\id})$
are mappings from $S^3$ to itself and $\tau\colon S\times S\to S\times S$ is the flip map, i.e., $\tau (x,y) =(y,x)$.
Throughout the paper we shall simply say that $(S,s)$ is a solution of the PE. We say that the solution $(S,s)$
is \emph{finite} if $S$ is a finite set, \emph{bijective} provided $s$ is a bijective map, \emph{bijective of finite order}
if there exists a positive integer $n$ such that $s^n=\id$, and \emph{involutive} if $s^2=\id$.

Writing $s(x,y)=(x\cdot y,\theta_x(y))$ for $x,y\in S$ it is easy to see that $(S,s)$ is a solution of the PE if and only if the following equalities
\begin{align}
	(x\cdot y)\cdot z & =x\cdot (y\cdot z),\label{eq:p1}\\
	\theta_x(y)\cdot\theta_{x\cdot y}(z) & =\theta_x(y\cdot z),\label{eq:p2}\\
	\theta_{\theta_x(y)}\theta_{x\cdot y} & =\theta_y\label{eq:p3}
\end{align}
hold for all $x,y,z\in S$. In particular, $(S,\cdotp)$ must be a semigroup (dealing with a solution $(S,s)$ of the PE we shall frequently
denote the multiplication in $S$ as a concatenation, that is $x\cdot y=xy$ for $x,y\in S$). Moreover, $(S,s)$ is a solution of the PE if
and only if the map $t=\tau s \tau\colon S\times S\to S\times S$ satisfies \[t_{12}t_{13}t_{23}=t_{23}t_{12},\]
where $t_{12}=t\times{\id}$, $t_{23}={\id}\times t$ and $t_{13}=(\tau\times{\id})({\id}\times t)(\tau\times{\id})$
are mappings from $S^3$ to itself. Such a map $t$ is called a \emph{solution of the Reversed Pentagon Equation (RPE)}.
Note that $(S,s)$ is a bijective solution of the PE if and only if $(S,s^{-1})$ is a bijective solution of the RPE. In particular,
an involutive solution of the PE also is a solution of the RPE.

\begin{exa}[{cf. \cite{CMM19x}}]
	Assume that $S$ is a semigroup and $f$ an idempotent endomorphism of $S$. Then the map $s\colon S\times S\to S\times S$,
	defined as $s(x,y)=(xy,f(y))$, is a solution of the PE. In particular, if $e\in S$ is an idempotent
	then $s(x,y)=(xy,e)$ is a solution of the PE.
\end{exa}

\begin{exa}[{cf. \cite{Mi98}}]
	If $S$ is a set and $f,g\colon S\to S$ are commuting idempotent maps then the map $s\colon S\times S\to S\times S$,
	given by $s(x,y )=(f(x),g(y))$, is a solution of both the PE and the RPE.
\end{exa}

\begin{exa}
	Let $G$ be a group of finite exponent. Let $E=\{1,\dotsc,n\}$ and let $\sigma\in\Sym(n)$ be a permutation satisfying 
	\begin{equation}\label{eq:sigma}
	\sigma^{\sigma(i)+1}=\sigma^i
	\end{equation}
	for each $i\in E$. Put $S=E\times G$ and let $s\colon S\times S\to S\times S$ be defined by \[s((i,a),(j,b))=((i,ab),(\sigma^i(j),b)).\]
	Then $(S,s)$ is a bijective solution of the PE of order equal to the least common multiple of the order of $\sigma$ and the exponent of $G$.
	First we shall check that $(S,s)$ is a solution of the PE. Observe that \eqref{eq:p1} is trivially satisfied because $S$
	is a semigroup for the multiplication defined by $(i,a)\cdot(j,b)=(i,ab)$. Moreover, 
	\[\theta_{(i,a)}(j,b)\theta_{(i,a)(j,b)}(k,c)= (\sigma^i(j),b)\theta_{(i,ab)}(k,c)=(\sigma^i(j),b)(\sigma^i(k),c)=(\sigma^i(j),bc)\]
	and \[\theta_{(i,a)}((j,b)(k,c))=\theta_{(i,a)}(j,bc)=(\sigma^i(j),bc),\]
	i.e., condition \eqref{eq:p2} holds. To prove that condition \eqref{eq:p3} holds, first note that $\sigma^{\sigma^i(j)+i}=\sigma^j$ for all 
	$i,j\in E$. Indeed, by assumption, the formula holds for $i=1$. So assume $i>1$. The result then follows by induction, because 
	\[\sigma^{\sigma^i(j)+i}=\sigma^{\sigma(\sigma^{i-1}(j))+1}\sigma^{i-1}=\sigma^{\sigma^{i-1}(j)}\sigma^{i-1}=\sigma^{\sigma^{i-1}(j)+i-1}.\]
	Hence, we also have that
	\begin{align*}
		\theta_{\theta_{(i,a)}(j,b)}\theta_{(i,a)(j,b)}(k,c) &= \theta_{(\sigma^i(j),b)}\theta_{(i,ab)}(k,c)
		=\theta_{(\sigma^i(j),b)}(\sigma^i(k),c)\\ & =(\sigma^{\sigma^i(j)+i}(k),c)=(\sigma^j(k),c)=\theta_{(j,b)}(k,c).
	\end{align*}
	Thus, indeed $(S,s)$ is a solution of PE. For a positive integer $m$ we get that 
	\[s^m((i,a),(j,b))=((i,ab^m),(\sigma^{mi}(j),b)).\]
	Hence $s^m=\id$ if and only if $m$ is a multiple of both the exponent of $G$ and the order of $\sigma$. Consequently, $s$ has finite
	order equal to the least common multiple of the order of $\sigma$ and the exponent of $G$. Finally, $(S,s)$ is also a solution of the
	RPE if and only if $G$ is the trivial group and $\sigma$ is either the identity or an involution. 
	
	A concrete example is the following: $E=\{1,2,3,4\}$ and $\sigma\in\{ \id, (1\,2)(3\,4), (1\,4)(2\,3), (1\,4\,3\,2)\}$.
	One can easily verify that actually on $E$ these $\sigma$ are (the unique) permutations in $\Sym(4)$ that satisfy
	condition \eqref{eq:sigma}. 
\end{exa}

The following natural definition is needed to determine when solutions are considered to be isomorphic
and it also yields a category of the solutions of the PE.

\begin{dfn}\label{def:hom}
	Let $(S,s)$ and $(S',s')$ be solutions of the PE. A morphism from $(S,s)$ to $(S',s')$ is a map $f\colon S\to S'$ such that
	$(f\times f)s=s'(f\times f)$ or, in other words, if the following diagram
	\[\begin{tikzcd}[row sep=30pt,column sep=30pt]
	S\times S\ar[r,"f{\times}f"]\ar[d,"s",swap] & S'\times S'\ar[d,"s'"]\\
	S\times S\ar[r,"f{\times}f"] & S'\times S'
	\end{tikzcd}\]
	is commutative. If a bijective morphism $f$ exists then we say that the solutions $(S,s)$ and $(S',s')$ are isomorphic.
\end{dfn}

With notations as in the Definition~\ref{def:hom}, write $s(x,y)=(x\cdot y,\theta_x(y))$ and $s'(a,b)=(a*b,\theta'_a(b))$.
One easily verifies that a map $f\colon S\to S'$ is a morphism of solutions if and only if
\begin{equation}\label{eq:hom}
f(x\cdot y)=f(x)*f(y)\qquad\text{and}\qquad f\theta_x=\theta'_{f(x)}f
\end{equation}
for all $x,y\in S$. In particular, $f$ is a morphism of semigroups.

Note also that if $(S_1,s_1)$ and $(S_2,s_2)$ are solutions of the PE then $(S,s)$, where $S=S_1\times S_2$ and the map
$s\colon S\times S\to S\times S$ is defined by the following formula \[s((x_1,x_2),(y_1,y_2))=(s_1(x_1,y_1),s_2(x_2,y_2)),\]
is again a solution of the PE. Clearly, $(S,s)$ is just a product of $(S_1,s_1)$ and $(S_2,s_2)$
in the category of solutions of the PE. Hence we may write $(S,s)=(S_1,s_1)\times(S_2,s_2)$.

\section{Bijective set-theoretic solutions of the Pentagon Equation}

In this section we show some basic properties of bijective solutions of finite order.
These will be crucial to deal with arbitrary involutive solutions.

It is clear that semigroups will play a role in the description of solutions of the PE.
For background and details on this topic we refer the reader to \cite{CP}.

Recall that a semigroup $S$ is called \emph{(left) simple} provided $S$ is the only (left) ideal of $S$.
Clearly each left simple semigroup is a simple semigroup.

\begin{lem}\label{lem:S=SS}
	Assume that $(S,s)$ is a bijective solution of the PE. Then $S=SS$ and $T=\{\theta_x:x\in S\}$
	is a subsemigroup of the transformation monoid $\Map(S,S)$.
	\begin{proof}
		Let $x,y \in S$. Since $s$ is bijective there exist $u,v \in S$ such that $s(u,v)=(x,y)$, i.e., $x=uv$ and 
		$y=\theta_u(v)$. Therefore, $S=SS$. Moreover, from \eqref{eq:p3} it follows that $\theta_y\theta_x=
		\theta_{\theta_u(v)}\theta_{uv}=\theta_v\in T$. Hence, $T$ is a semigroup.
	\end{proof}
\end{lem}

\begin{pro}\label{pro:sgr}
	Assume that $(S,s)$ is a solution of the PE. If $s$ is bijective of finite order then $S$ is a simple semigroup.
	In particular, $|S|=1$ if and only if $S$ has a zero element.
	\begin{proof}
		Since $s$ is of finite order there exists $n \in \mathbb{N}$ such that $s^n=\id$. Then, for each $x,y\in S$,
		there exists $z\in S$ such that $xyz=x$. Hence, for any $x,y\in S$, we have $x\in SyS$. Therefore,
		every principal ideal, and thus every ideal of $S$, equals $S$, i.e., $S$ is a simple semigroup. 
	\end{proof}
\end{pro}

Let $S$ be a semigroup and $\theta_x\colon S\to S$ a map for every $x\in S$. Consider the map $s\colon S\times S \to S \times S$
defined by $s(x,y)=(xy,\theta_x(y))$. Then $s$ is a solution of the RPE if and only if the following equalities 
\begin{align}
	xy & =x\theta_y(z)yz,\label{eq:rp1}\\
	\theta_{x\theta_y(z)}(yz) & =\theta_x(y)z,\label{eq:rp2}\\
	\theta_{\theta_x(y)} & =\theta_x\theta_y\label{eq:rp3}
\end{align}
hold for all $x,y,z\in S$.

Recall that a semigroup $E$ is called a \emph{left zero semigroup} if every element of $E$ is a left zero of $E$, i.e., $xy=x$ for all $x,y\in E$.
Moreover, a semigroup $S$ is said to be a \emph{left group} if $S$ is a direct product $E\times G$ of a left zero semigroup $E$ and a group $G$.
Obviously, a left group is left simple, and it is easy to see that each left group is right cancellative, whereas each idempotent of a left group
is a right identity (see \cite{CP} for more details). 

\begin{pro}\label{pro:GxE}
	Let $(S,s)$ be a solution of both the PE and the RPE (for example this is the case if $(S,s)$ is an involutive solution of the PE).
	If $s$ is bijective of finite order then $S$ is a left group.
	\begin{proof}
		Since $s$ is of finite order, there exists a positive integer $n$ such that $s^n=\id$. 
		Let $x \in S$. Clearly, there exists $e\in S$ such that $xe=x$. By Lemma~\ref{lem:S=SS}, there exist $u,v\in S$
		such that $x=uv$. Hence, since $(S,s)$ also is a solution of the RPE, from \eqref{eq:rp1} we get that
		$x=uv=u\theta_v(z)vz$ for every $z\in S$. Hence $S=Sz$ for each $z\in S$. So, $S$ is left simple.
		Therefore, $e\in S=Sx$ and there exists $z\in S$ such that $e=zx$. So, $e^2=(zx)e=z(xe)=zx=e$, that is, $e$ is an idempotent.
		Thus $S$ is left simple and contains an idempotent. Consequently, by \cite[Theorem I.1.27]{CP}, we conclude that $S$ is a left group.
	\end{proof}
\end{pro}

\begin{lem}\label{lem:left-0-th}
	Let $(S,s)$ be a bijective solution of the PE. If $S$ is a left zero semigroup then all maps $\theta_x$ for $x\in S$ are bijective.
	Moreover, for each $x\in S$, the map $\theta_x$ is either the identity or a fixed-point free permutation.
	\begin{proof}
		Since $S$ is a left zero semigroup, we have $s(x,y)=(x,\theta_x(y))$ for all $x,y \in S$. If $\theta_x(y)=\theta_x(z)$ for some $x,y,z\in S$
		then $s(x,y)=(x,\theta_x(y))=(x,\theta_x(z))=s(x,z)$. Since $s$ is bijective, we get $y=z$ and thus $\theta_x$ is injective.
		Moreover, let $u\in S$ be such that $s(x,u)=(x,y)$. Then $\theta_x(u)=y$ and thus $\theta_x$ also is surjective. 
		
		Furthermore, if $\theta_x$ has a fixed point, say $z\in S$, then, by \eqref{eq:p3}, 
		$\theta_z\theta_x=\theta_{\theta_x(z)}\theta_{xz}=\theta_z$. Because $\theta_z$ is bijective this leads to $\theta_x=\id$, as desired.
	\end{proof}
\end{lem}

\section{Involutive solutions: a reduction}

In this section we will prove that the description of all involutive set-theoretic solutions $(S,s)$ of the PE on a semigroup $S$ 
can be reduced to the description of solutions on a left zero semigroup. More precisely, we will show that $S$ may be identified with
$E\times G$, where $E$ is a left zero semigroup and $G$ is a group, and then that $s$ is composed of solutions $s_E$ and $s_G$
on $E$ and $G$, respectively. The solution $s_G$ is unique as it is known that on a group there only is one bijective (not necessarily
involutive) solution of the PE. This has been shown by Kashaev and Sergeev in \cite{KaSe98}. Catino, Mazzotta and Micolli extended this
result by showing that all solutions of PE on a group are determined by normal subgroups \cite{CMM19x}. For completeness' sake we give
an easy proof that on a group there is unique bijective solution of the PE. 

\begin{pro}[{see \cite{KaSe98,CMM19x}}]\label{GroupUnique}
	Assume that $(G,s)$ is a bijective solution of the PE, where $G$ is a group. Then $s(x,y)=(xy,y)$ for all $x,y\in G$.
	In particular, there is a unique solution of the PE on a group $G$, which will be denoted $s_G$. Moreover, the order of
	$s$ equals the exponent of $G$; and thus $s$ is involutive precisely when $G$ is an elementary abelian $2$-group.
	\begin{proof}
		As before write $s(x,y)=(xy,\theta_x (y))$ for $x,y\in G$. Let $1$ be the identity element of the group $G$.
		First we show that $\theta_1$ is injective. Indeed, for $x\in G$ note that, by \eqref{eq:p2},
		\[\theta_x (x^{-1}) =\theta_x (x^{-1}1)=\theta_x(x^{-1}) \theta_{xx^{-1}} (1) =\theta_x(x^{-1}) \theta_1(1).\]
		Thus, since $G$ is a group, it follows that $\theta_1(1)=1$. Hence, again by \eqref{eq:p2},
		\[1=\theta_1 (1)=\theta_{1}(xx^{-1})=\theta_{1}(x) \theta_x(x^{-1}).\]
		Consequently, $\theta_{1}(x)^{-1}=\theta_x(x^{-1})$.
		
		Now, assume that $\theta_{1}(x) =\theta_{1}(y)$ for some $x,y\in G$. Then, by the previous,
		\[\theta_x(x^{-1})=\theta_{1}(x)^{-1} =\theta_{1}(y)^{-1}=\theta_y(y^{-1}).\]
		Therefore, \[s(x,x^{-1}) =(1,\theta_x(x^{-1}))=(1, \theta_y(y^{-1})) =s(y,y^{-1}).\]
		Since, by assumption, $s$ is bijective, it follows that $x=y$ and thus, indeed, $\theta_{1}$ in injective.
		
		From \eqref{eq:p2} we also obtain that $\theta_x(1) \theta_x(1) =\theta_x(1)$ for all $x\in G$ and thus
		$\theta_x (1)=1$ for all $x\in G$. Hence, \eqref{eq:p3} yields that $\theta_1 \theta_x=\theta_{\theta_x(1)}\theta_x =\theta_1$.
		Since, by the above, $\theta_1$ is injective, we obtain that $\theta_x =\id$ for all $x\in G$. Consequently, $s(x,y)=(xy,y)$ for
		all $x,y\in G$.
	\end{proof}
\end{pro}

Recall that a solution $(S,s)$ of the PE is called involutive provided $s^2=\id$ or, equivalently, the following equalities
\begin{align}
	xy\theta_x(y) & =x,\label{eq:i1}\\
	\theta_{xy}(\theta_x(y)) & =y\label{eq:i2}
\end{align}
hold for all $x,y\in S$.

\begin{thm}\label{thm:involuiveGxI}
	Let $(S,s)$ be an involutive solution of the PE. Then there exist a left zero semigroup $E$, an involutive solution $(E,s_E)$ of the PE
	on $E$, and an elementary abelian $2$-group $G$ such that $S$ may be identified with $E\times G$ and then \[(S,s)=(E,s_E) \times (G,s_G),\]
	where $s_G$ is the unique bijective solution of the PE on the group $G$. Moreover, $T=\{\theta_x:x \in S\}$ is an elementary abelian $2$-group.
	The converse is also true, that is, if $(E,s_E)$ is an involutive solution of the PE on a left zero semigroup $E$ and $(G,s_G)$ is
	the unique bijective solution of the PE on an elementary abelian $2$-group $G$ then the product $(E,s_E)\times(G,s_G)$ is an involutive
	solution of the PE.
	\begin{proof} 
		Note that $(S,s)$ also is a solution of the RPE and thus, by Proposition~\ref{pro:GxE}, $S$ may be identified with $E\times G$,
		where $E$ is a left zero semigroup and $G$ is a group. Moreover, $E(S)=E\times \{1\}$, where by $E(S)$ we denote the set of
		idempotents of the semigroup $S$. Recall also that each idempotent of $S$ is a right identity in $S$, that is if $e\in E(S)$
		then $xe=x$ for each $x\in S$. Let $x \in S$ and $e\in E(S)$. Then, by \eqref{eq:p2}, we obtain 
		\[\theta_x(e)=\theta_x(e^2)=\theta_x(e)\theta_{xe}(e)=\theta_x(e)\theta_x(e),\]
		i.e., $\theta_x(e)\in E(S)$. If $e=(i,1)$ for some $i \in E$ then we denote by $\vartheta_x(i)$ the unique element of $E$
		such that $\theta_x(e)=(\vartheta_x(i),1)$. If $y=(i,g)\in S$ then \eqref{eq:p2} yields
		\[\theta_x(y)=\theta_x(ey)=\theta_x(e)\theta_{xe}(y)=\theta_x(e)\theta_x(y)=(\vartheta_x(i), 1) \theta_x(y)\]
		and thus $\theta_x(y)=(\vartheta_x(i),h)$ for some $h\in G$. Moreover, write $x=(j,a)$ for some $j\in E$ and $a\in G$.
		Then it follows, by \eqref{eq:i1}, that $(j, agh)=xy\theta_x(y)=x =(j,a)$. Hence we have $h=g^{-1}$ and consequently 
		$\theta_x(i,g)=\theta_x(y)=(\vartheta_x(i), g^{-1})$. Finally, by \eqref{eq:p3},
		\[(\vartheta_{(\vartheta_x(i),1)}\vartheta_x(j),g)
		=\theta_{\theta_x(e)}\theta_x(y)=\theta_{\theta_x(e)}\theta_{xe}(y)=\theta_e(y)=(\vartheta_e(i),g^{-1}).\]
		So, $g=g^{-1}$ for all $g\in G$, and thus $G$ is an elementary abelian $2$-group. We obtained
		\[\theta_x(i,g)=(\vartheta_x(i),g).\] Moreover, note that, by \eqref{eq:i2},
		$(\vartheta^2_x(i),1)=\theta^2_x(i,1)=\theta_{x(i,1)}\theta_x(i,1) =(i,1)$, i.e., $\vartheta_x^2(i)=i$. Hence, we obtain
		$\theta_x^2(i,g)=(\vartheta^2_x(i),g)=(i,g)$, i.e., $\theta^2_x=\id$. It follows that $T=\{\theta_x:x \in S\}$ is an elementary abelian $
		2$-group, since $T$ is a semigroup by Lemma~\ref{lem:S=SS} and, by \eqref{eq:rp3}, $\id =\theta^2_x=\theta_{\theta_x(x)}\in T$.
		
		Next, we claim that the map $\vartheta_x$ does not depend on $a$ (recall that $x=(j,a)$), that is, $\vartheta_{(j,a)}=\vartheta_{(j,1)}$.
		To prove this, first note that $\theta_{\theta_x(y)}=\theta_x\theta_y$ by \eqref{eq:rp3} and also $\theta_{\theta_x(y)}\theta_{xy}=\theta_y$
		by \eqref{eq:p3}. Therefore, we get $\theta_{xy}=(\theta_x\theta_y)^{-1}\theta_y=\theta_x$, since $T$ is an elementary abelian $2$-group.
		In particular, \[ (\vartheta_x(i),1)=\theta_x(i,1)=\theta_{(j,1)x}(i,1) =\theta_{(j,1)}(i,1)=(\vartheta_{(j,1)}(i),1),\]
		i.e., $\vartheta_{(j,a)}(i)=\vartheta_{(j,1)}(i)$. This proves the claim.
		
		For the sake of simplicity we may identify the elements of $E(S)$ with the elements of $E$ and write $\vartheta_{(j,1)}$ simply as 
		$\vartheta_j$. Therefore, \[\theta_{(i,a)}(j,g)=\theta_{(i,1)}(j,g)=(\vartheta_i(j),g).\]
		Now consider $s_E$, the restriction of $s$ to $E\times E$. We claim that $(E,s_E)$ is an involutive solution of the PE. First,
		note that $s_E(i,j)=(i,\vartheta_i(j))$ for $i,j\in E$. Hence, \eqref{eq:p1} follows since $E$ is a subsemigroup of $S$,
		\eqref{eq:p2} trivially follows since $E$ is a left zero semigroup (i.e., $ij=i$ for all $i,j \in E$). Let $i,j,k \in E$. Then
		\[(\vartheta_{\vartheta_i(j)}\vartheta_i(k),1) =( \vartheta_{\vartheta_i(j)}\vartheta_{ij}(k),1)= 
		\theta_{\theta_{(i,1)}(j,1)}\theta_{(i,1)(j,1)}(k,1)=\theta_{(j,1)}(k,1)=( \vartheta_j(k),1),\]
		which means that $\vartheta_{\vartheta_i(j)}\vartheta_i=\vartheta_j$, that is \eqref{eq:p3} holds. Moreover, since $\vartheta_i^2=\id$,
		it holds that $s_E^2(i,j)=s_E(i,\vartheta_i(j))=(i,\vartheta_i^2(j))=(i,j)$. Therefore, $(E,s_E)$ is an involutive solution, as claimed.
		
		Finally, from Proposition~\ref{GroupUnique} we know that on $G$ there is a unique bijective solution defined by $s_G(g,h)=(gh,h)$
		for $g,h\in G$. Since $G$ is an elementary abelian $2$-group, $s_G$ is involutive. Then, 
		\[s((i,g),(j,h))=((i,g)(j,h),\theta_{(i,g)}(j,h))=((i,gh),(\vartheta_i(j),h))=(s_E\times s_G)((i,j),(g,h)).\]
		Since the converse readily can be verified, the result follows.
	\end{proof}
\end{thm}

It follows from the previous result that if $(S,s)$ is an involutive solution then $S= E\times G$ and $G$ is an elementary
abelian $2$-group. Hence, $E\times\{1\}=E(S)=\{x^2:x\in S\}$ and $G$ is isomorphic to $x^2Sx^2$ for any $x \in S$.

\begin{cor}\label{cor:basicFact}
	Let $(S,s)$ be an involutive solution of the PE. Then:
	\begin{enumerate}
		\item $\theta_{\theta_x(y)}=\theta_x\theta_y$, $\theta^2_x=\id$ and $\theta_{xy}=\theta_x$ for all $x,y\in S$.
		\item $x\theta_y(z)=xz$ for all $x,y,z\in S$.
		\item $\theta_x$ is an automorphism of $S$ for each $x\in S$. 
	\end{enumerate}
	\begin{proof}
		From Theorem~\ref{thm:involuiveGxI} we know that $S=E\times G$ with $G$ an elementary abelian $2$-group and 
		$E= E(S)=\{x^2:x\in S\}$. Moreover, we have $x^3=x$ for each $x\in S$. Furthermore, from the proof of
		Theorem~\ref{thm:involuiveGxI}, we know that $\theta_x=\theta_{x^2}$ and $\theta_x^{2}=\id$.
		
		Since $(S,s)$ is an involutive solution of the PE and thus also of the RPE, it follows from from \eqref{eq:rp3} that 
		$\theta_{\theta_x(y)}=\theta_x\theta_y$ for all $x,y\in S$. From \eqref{eq:p3} we thus also get that
		\[\theta_{xy}=\theta_{\theta_x(y)}^{-1}\theta_y=\theta_{\theta_x(y)}\theta_y=\theta_x\theta_y\theta_y=\theta_x\] for all $x,y \in S$. 
		This proves part (1).
		
		To prove part (2), note that Theorem~\ref{thm:involuiveGxI} yields that $y^2$ and $z^2$ are idempotents that are right identities in $S$. 
		It then follows from $\theta_y=\theta_{y^2}$ and \eqref{eq:rp1} that
		\[x\theta_y(z)=x\theta_{y^2}(z)y^2z^2=(x\theta_{y^2}(z)y^2z)z=xy^2z=xz.\]
		
		Part (3) follows from \eqref{eq:p2} and part (1). Indeed, $\theta_x(yz)=\theta_x(y)\theta_{xy}(z)=\theta_x(y)\theta_x(z)$.
	\end{proof}
\end{cor}

Since the so-called commutative and cocommutative solutions of the PE play a role, especially in the context
of locally compact groups and Hilbert spaces (see \cite{BaSk93}), we recall definitions of these notions and give
a simple proof that involutive solutions are of this type. A solution $(S,s)$ of the PE is said to be \emph{commutative}
if $s_{12}s_{13}=s_{13}s_{12}$, i.e., if the following equalities
\begin{align}
	xyz&=xzy,\label{eq:c1}\\
	\theta_x&=\theta_{xy}\label{eq:c2}
\end{align}
hold for all $x,y,z \in S$. Similarly, $(S,s)$ is said to be \emph{cocommutative} provided
$s_{13}s_{23}=s_{23}s_{13}$, i.e., if the following equalities
\begin{align}
	x\theta_y(z)&=xz,\label{eq:cc1}\\
	\theta_x\theta_y&=\theta_y\theta_x\label{eq:cc2}
\end{align}
hold for all $x,y,z \in S$. 

\begin{cor}
	Assume that $(S,s)$ is an involutive solution of the PE. Then $(S,s)$ is both commutative and cocommutative. 
	\begin{proof}
		By Theorem~\ref{thm:involuiveGxI} we may assume that $S=E\times G$ is the direct product of a left zero semigroup $E$
		and an elementary abelian $2$-group $G$. Let $(e_i,g_i)\in S$ for $1\leq i\leq 3$. Then
		\[(e_1,g_1)(e_2,g_2)(e_3,g_3)=(e_1,g_1g_2g_3)=(e_1,g_1g_3g_2)=(e_1,g_1)(e_3,g_3)(e_2,g_2),\]
		i.e., \eqref{eq:c1} holds in $S$. Moreover, by Corollary~\ref{cor:basicFact}(1), \eqref{eq:c2} also holds.
		Hence $(S,s)$ is commutative. Furthermore, by Corollary~\ref{cor:basicFact}(2), we have that \eqref{eq:cc1} holds.
		Finally, \eqref{eq:cc2} holds since, by Theorem~\ref{thm:involuiveGxI}, $T=\{\theta_x:x\in S\}$ is an abelian group.
		Therefore, $(S,s)$ is cocommutative.
	\end{proof}
\end{cor}

\section{Retraction of an involutive solution and irretractable solutions}\label{subsec:ret}

In this section we introduce the notion of an irretractable involutive solution of the PE and we prove that
for a given cardinality there is unique such solution and we give an actual construction of this solution.

Assume that $(S,s)$ is an involutive solution of the PE. We define the equivalence relation $\sim$ on $S$,
called the {\em retract}, as follows. For $x,y\in S$, \[x\sim y\iff\theta_x=\theta_y.\]

Put $\overline{S}=S/{\sim}$ and let $\overline{x}$ denote the $\sim$-class of $x\in S$, called the retract class of $x$. 
From Corollary~\ref{cor:basicFact}(1) we know that $\theta_{xy}= \theta_x$ for all $x,y\in S$ or, in other words,
\begin{equation}\label{eq:star}
xy\sim x.
\end{equation}
In particular, if $x_1\sim x_2$ and $y_1\sim y_2$ for some elements $x_1,x_2,y_1,y_2\in S$ then $x_1y_1\sim x_1\sim x_2 \sim x_2y_2$.
Therefore, $\sim$ is a congruence on $S$ and $\overline{S}$ is a semigroup. Moreover, again by Corollary~\ref{cor:basicFact}(1),
we get \[\theta_{\theta_{x_1}(y_1)}=\theta_{x_1}\theta_{y_1}=\theta_{x_2}\theta_{y_2}=\theta_{\theta_{x_2}(y_2)}.\]
Hence, $\theta_{x_1}(y_1)\sim\theta_{x_2}(y_2)$ and thus the map $\overline{\theta}_{\overline{x}}\colon\overline{S}\to\overline{S}$,
given by $\overline{\theta}_{\overline{x}}(\overline{y})=\overline{\theta_x(y)}$, is well defined. 
This allows us to define a map $\overline{s}\colon\overline{S}\times \overline{S}\to \overline{S}\times \overline{S}$ by the formula
\[\overline{s}(\overline{x},\overline{y})=(\overline{x},\overline{\theta}_{\overline{x}}(\overline{y})).\]
Note that \eqref{eq:star} assures that $\overline{S}$ is a left zero semigroup. 
Hence the \emph{retract} \[\Ret(S,s)=(\overline{S},\overline{s})\] of $(S,s)$ satisfies both \eqref{eq:p1} and \eqref{eq:p2}. Further,
\[\overline{\theta}_{\overline{\theta}_{\overline{x}}(\overline{y})}\overline{\theta}_{\overline{x}}(\overline{z})
=\overline{\theta}_{\overline{\theta_x(y)}}(\overline{\theta_x(z)})=\overline{\theta_{\theta_x(y)}\theta_x(z)}=\overline{\theta_y(z)}
=\overline{\theta}_{\overline{y}}(\overline{z})\]
for all $x,y,z\in S$, i.e., the equality \eqref{eq:p3} also holds for $\Ret(S,s)$. So, $\Ret (S,s)$ is a solution of the PE.
Finally, by Corollary~\ref{cor:basicFact}(1), $\theta_x^2=\id$ for each $x\in S$; hence
also $\overline{\theta}_{\overline{x}}^2=\id$. It follows that $\Ret(S,s)$ satisfies both equalities \eqref{eq:i1}
and \eqref{eq:i2}. Thus we have proved the first part of Proposition~\ref{pro:RII}.
For the second part we need the following definition.

\begin{dfn}
	An involutive solution $(S,s)$ of the PE is said to be \emph{irretractable} if $(S,s)=\Ret(S,s)$.
\end{dfn}

It is clear that a finite involutive solution $(S,s)$ of the PE is irretractable if and only if $|S|=|T|$, where
$T=\{\theta_x:x \in S\}$. In this case, by Theorem~\ref{thm:involuiveGxI}, $|S|$ is a power of $2$.

\begin{pro}\label{pro:RII}
	Let $(S,s)$ be an an involutive solution of the PE. Then $\Ret(S,s)$ is an irretractable involutive solution of the PE.
	\begin{proof}
		We only need to show the second part. This is shown as follows. If $x,y\in S$ then by Corollary~\ref{cor:basicFact}(1) we get
		\begin{align*}
			\overline{\theta}_{\overline{x}}=\overline{\theta}_{\overline{y}}
			& \iff \overline{\theta_x(z)}=\overline{\theta_y(z)}\text{ for all }z\in S\\
			& \iff \theta_{\theta_x(z)}=\theta_{\theta_y(z)}\text{ for all }z\in S\\
			& \iff \theta_x\theta_z=\theta_y\theta_z\text{ for all }z\in S\\
			& \iff \theta_x=\theta_y\\
			& \iff \overline{x} =\overline{y}.
		\end{align*}
		Hence the result is proved.
	\end{proof}
\end{pro}

\begin{lem}\label{lem:irr-l-0}
	Assume that $(S,s)$ is an involutive solution of the PE. Then $\theta_x(y)\sim\theta_y(x)$ for all $x,y\in S$. If in addition
	the solution $(S,s)$ is irretractable then $\theta_x(y)=\theta_y(x)$ for all $x,y\in S$ and $S$ is a left zero semigroup.
	\begin{proof}
		Let $x,y \in S$. By Theorem~\ref{thm:involuiveGxI} we know that $T=\{\theta_x:x\in S\}$ is an abelian group. Then,
		by Corollary~\ref{cor:basicFact}(1), we get $\theta_{\theta_x(y)}=\theta_x\theta_y=\theta_y\theta_x=\theta_{\theta_y(x)}$.
		Hence $\theta_x(y) \sim \theta_y(x)$. This proves the first part of the statement.
		
		For the second part, assume that $(S,s)$ also is irretractable. By Theorem~\ref{thm:involuiveGxI} we may assume that
		$S= E\times G$, where $E$ is a left zero semigroup and $G$ is an elementary abelian $2$-group. Let $e\in E$ and $g\in G$.
		Then Corollary~\ref{cor:basicFact}(1) yields $\theta_{(e,g)}=\theta_{(e,1)(e,g)}=\theta_{(e,1)}$.
		Hence, the irretractability of $(S,s)$ gives $g=1$. Thus $|G|=1$ and $S= E\times\{1\}\cong E$ is a left zero semigroup
	\end{proof}
\end{lem}

In next two propositions we completely describe the structure of all irretractable involutive solutions of the PE.

\begin{pro}\label{pro:irretract-2-abelian}
	Let $(A,+)$ be an elementary abelian $2$-group. Define the map $t\colon A\times A\to A\times A$ by $t(x,y)=(x,x+y)$.
	Then $(A,t)$ is an irretractable involutive solution of the PE. Conversely, if $(S,s)$ an irretractable involutive 
	solution of the PE then there exists a natural structure $(S,+)$ of an elementary abelian $2$-group on $S$ such that
	$s(x,y)=(x,x+y)$ for all $x,y\in S$.
	\begin{proof}
		First, note that we have $t=\tau s_A \tau$, where $s_A$ is the unique solution of the PE on the group $A$
		(see Proposition~\ref{GroupUnique}) and $\tau\colon A\times A\to A\times A$ is the flip map. Hence, $t$
		is an involutive solution of the RPE and the PE. It remains to prove that $(A,t)$ is irretractable. If $x,y\in A$ then
		\begin{align*}
			\theta_x=\theta_y & \iff \theta_x(z)=\theta_y(z)\text{ for all }z\in A\\
			& \iff x+z=y+z\text{ for all }z\in A\\
			&\iff x=y.
		\end{align*}
		Therefore, $\Ret(A,t)=(A,t)$ and the solution $(A,t)$ is irretractable.
		
		Conversely, assume that $(S,s)$ is an irretractable involutive solution of the PE.
		By Lemma~\ref{lem:irr-l-0} we know that $s$ is of the form $s(x,y)=(x,\theta_x(y))$.
		We claim that the operation $x+y=\theta_x(y)$ for $x,y\in S$ makes $(S,+)$ an elementary
		abelian $2$-group. Observe first that, by Corollary~\ref{cor:basicFact}(1), we get
		\[x+(y+z)=x+\theta_y(z)=\theta_x\theta_y(z)=\theta_{\theta_x(y)}(z)=\theta_{x+y}(z)=(x+y)+z.\]
		Hence $(S,+)$ is a semigroup. Further, by Lemma~\ref{lem:irr-l-0}, \[x+y=\theta_x(y)=\theta_y(x)=y+x,\]
		i.e., the semigroup $(S,+)$ is commutative. Since the solution $(S,s)$ is irretractable, there exists a unique element $0\in S$
		such that $\theta_0=\id$ (it is enough to put $0=\theta_a(a)$ for any $a\in S$ because 
		$\theta_{\theta_a(a)}=\theta_a^2=\id$). Then \[0+x=\theta_0(x)=x\qquad\text{and}\qquad x+0=0+x=x.\] 
		Hence $0$ is the identity of $(S,+)$.
		Finally, $x+x=\theta_x(x)=0$. Therefore, $(S,+)$ is a elementary abelian $2$-group, as claimed, and $s(x,y)=(x,x+y)$.
	\end{proof}
\end{pro}

The solution $(A,t)$ on an elementary abelian $2$-group $(A,+)$ will be denoted as $(A,t_A)$.
Recall also that $t_A=\tau s_A\tau$, where $\tau\colon A\times A\to A\times A$ is the flip map
and $s_A$ is the unique bijective solution of the PE defined on the group $(A,+)$.

Next we show that irretractable involutive solutions of the PE are uniquely, up to isomorphism, determined by their cardinalities.

\begin{pro}\label{dim-arg}
	Assume that $(S,s)$ and $(S',s')$ are irretractable involutive solutions of the PE.
	Then the solutions $(S,s)$ and $(S',s')$ are isomorphic if and only if $|S|=|S'|$.
	\begin{proof}
		Clearly if $(S,s)$ and $(S',s')$ are isomorphic then $|S|=|S'|$. Conversely, assume that $|S|=|S'|$.
		Since $(S,s)$ and $(S',s')$ are irretractable involutive solutions, Proposition~\ref{pro:irretract-2-abelian} implies
		that $S$ and $S'$ admit natural structures of elementary abelian $2$-groups, let us denote them by $(S,+)$
		and $(S',+)$ respectively, such that $s(x,y)=(x,x+y)$ for all $x,y\in S$ and $s'(a,b)=(a,a+b)$ for all $a,b\in S'$.
		We claim that the groups $(S,+)$ and $(S',+)$ are isomorphic. This is obvious in case $|S|=|S'|$ is finite, as then,
		$|S|=|S'|=2^n$ and $(S,+)\cong C_2^n\cong(S',+)$ for some non-negative integer $n$. So, assume that $S$ and $S'$ are infinite.
		Considering $S$ and $S'$ as vector spaces over the field $\mathbb{F}_2$ with two elements, we get 
		\[\dim_{\mathbb{F}_2}S=|S|=|S'|=\dim_{\mathbb{F}_2}S'.\]
		Hence, it follows that groups $(S,+)$ and $(S',+)$ are isomorphic, and the claim is proved.
		
		Therefore, there exists a group isomorphism $f\colon(S,+)\to(S',+)$. Moreover,
		\begin{align*}
			(f\times f)s(x,y) & =(f\times f)(x,x+y)=(f(x),f(x+y))\\
			& =(f(x),f(x)+f(y))=s'(f(x),f(y))=s'(f\times f)(x,y),
		\end{align*}
		which means that $f$ is also a morphism of solutions. This finishes the proof.
	\end{proof}
\end{pro}

\section{A description of all involutive solutions of the PE}

In this section we give a description of all involutive solutions of the PE. We shall start with the following construction.

\begin{pro}\label{pro:retrsol-prime}
	Let $(A,t_A)$ be the irretractable involutive solution of the PE defined on an elementary abelian $2$-group $(A,+)$.
	Assume that $X$ is a non-empty set and $\sigma\colon A\to\Sym(X)$. Put $S=X\times A$ and define the map
	$s\colon S\times S\to S\times S$ by the formula \[s((x,a),(y,b)) =((x,a),(\sigma_{a+b}\sigma_b^{-1}(y),a+b)).\]
	Then $(S,s)$ is an involutive solution of the PE. Moreover, $\Ret(S,s)=(A,t_A)$.
	\begin{proof}
		Write $\theta_{(x,a)}(y,b)=(\sigma_{a+b}\sigma_b^{-1}(y),a+b)$. Define $(x,a)\cdot(y,b)=(x,a)$ for $(x,a)\in S$ and $(y,b)\in S$.
		Then $(S,\cdotp)$ is a left zero semigroup. Hence both equalities \eqref{eq:p1} and \eqref{eq:p2} hold for $(S,s)$. Moreover, we have
		\begin{align*}
			\theta_{\theta_{(x,a)}(y,b)}\theta_{(x,a)\cdot(y,b)}(z,c)
			& =\theta_{(\sigma_{a+b}\sigma_b^{-1}(y),a+b)}\theta_{(x,a)}(z,c)\\
			& =\theta_{(\sigma_{a+b}\sigma_b^{-1}(y),a+b)}(\sigma_{a+c}\sigma_c^{-1}(z),a+c)\\
			& =(\sigma_{(a+b)+(a+c)}\sigma_{a+c}^{-1}\sigma_{a+c}\sigma_c^{-1}(z),(a+b)+(a+c))\\
			& =(\sigma_{b+c}\sigma_c^{-1}(z),b+c)\\
			& =\theta_{(y,b)}(z,c).
		\end{align*}
		Thus, also the equality \eqref{eq:p3} holds for $(S,s)$. Furthermore,
		\begin{align*}
			s^2((x,a),(y,b)) &=s((x,a),(\sigma_{a+b}\sigma_b^{-1}(y),a+b))\\
			& =((x,a),(\sigma_{a+(a+b)}\sigma_{a+b}^{-1}\sigma_{a+b}\sigma_b^{-1}(y),a+(a+b)))\\
			& =((x,a),(y,b)).
		\end{align*}
		Therefore, $(S,s)$ is an involutive solution of the PE.
		
		Finally, if $(x,a)\in S$ and $(y,b)\in S$ then $(x,a)\sim(y,b)$ or, in other words, $\theta_{(x,a)}=\theta_{(y,b)}$
		if and only if $a=b$. Hence the quotient $S/{\sim}$ may be identified with $A$. It remains to observe that under this
		identification we get $\Ret(S,s)=(A,t_A)$. Thus the result is proved.
	\end{proof}
\end{pro}

The solution $(S,s)$ constructed in Proposition~\ref{pro:retrsol-prime} will be called the
\emph{extension of $(A,t_A)$ by $X$ and $\sigma$} and denoted by \[\Ext_X^\sigma(A,t_A).\]

\begin{lem}\label{lem:classes}
	Let $(S,s)$ be an involutive solution of the PE. Then all retract classes on $S$ have the same cardinality.
	\begin{proof}
		Fix $x,y\in S$. Define $X=\{z\in S:\theta_z=\theta_x\}$ and $Y=\{z\in S:\theta_z=\theta_y\}$, the retract classes of
		$x$ and $y$, respectively. Consider the map $f\colon X\to Y$ given by the formula $f(z)=\theta_y\theta_x(z)$.
		Note that $f$ is well defined, since for $z\in X$ we have, by Corollary~\ref{cor:basicFact}(1) and Theorem~\ref{thm:involuiveGxI}, that  
		\[\theta_{f(z)}=\theta_{\theta_y\theta_x(z)}=\theta_y\theta_{\theta_x(z)}= \theta_y\theta_x\theta_z=\theta_y\theta_x^2=\theta_y,\]
		and thus $f(z)\in Y$. Since $\theta_y\theta_x$ is a bijective map, clearly $f$ is injective.
		Moreover, if $w\in Y$ then put $z=\theta_x\theta_y(w)$. Then, again by Corollary~\ref{cor:basicFact}(1) and Theorem~\ref{thm:involuiveGxI},
		\[\theta_z= \theta_{\theta_x\theta_y(w)} =\theta_x\theta_{\theta_y(w)}=\theta_x\theta_y\theta_w=\theta_x\theta_y^2=\theta_x.\]
		So $z \in X$. Finally, $f(z)=\theta_y\theta_x\theta_x\theta_y(w)=\theta_x^2\theta_y^2(w)=w$. Thus $f$ is also surjective.
	\end{proof}
\end{lem}

The following result shows that solutions constructed in Proposition~\ref{pro:retrsol-prime}
are actually all solutions of the PE defined on left zero semigroups.

\begin{pro}\label{pro:semidir-prime}
	Assume that $(S,s)$ is an involutive solution of the PE defined on a left zero semigroup $S$.
	Then there exists an abelian elementary $2$-group $(A,+)$ and a non-empty set $X$ such that
	$S$ may be identified with $X\times A$ and then $(S,s)=\Ext_X^\sigma(A,t_A)$ for some $\sigma\colon A\to\Sym(X)$.
	\begin{proof}
		Write $s(x,y)=(x,\theta_x(y))$. By Proposition~\ref{pro:irretract-2-abelian} we know that if
		$\Ret(S,s)=(A,\overline{s})$ then $A$ is an elementary abelian $2$-group such that
		\[\overline{s}(a,b)=(a,\overline{\theta}_a(b))=(a,a+b)\] for all $a,b\in A$. Observe also that $(A,\overline{s})=(A,t_A)$.
		
		Next, choose $w\in S$ and define $X=\{x\in S:\theta_x=\theta_w\}$, the retract class of $w$.
		We know from Lemma~\ref{lem:classes} that there exist a bijection between $X$ and each retract
		class $\overline{x}$ of $x\in S$. Hence, we may identify $S$ with $X\times A$, and under this
		identification elements $(x,a)\in S$ and $(y,b)\in S$ lie in the same retract class if and only if
		$a=b$. Furthermore, it is clear that we may write \[s((x,a),(y,b))=((x,a),(\pi_{x,a,b}(y),a+b)),\]
		i.e., $\theta_{(x,a)}(y,b)=(\pi_{x,a,b}(y),a+b)$ for some maps $\pi_{x,a,b}\in\Map(X,X)$. 
		We claim that $\pi_{x,a,b}$ does not depend on $x$. To prove this choose $y\in X$.
		Since $(x,a)$ and $(y,a)$ lie in the same retract class, we get $\theta_{(x,a)}=\theta_{(y,a)}$. Hence
		\[(\pi_{x,a,b}(z),a+b)=\theta_{(x,a)}(z,b)=\theta_{(y,a)}(z,b)=(\pi_{y,a,b}(z),a+b)\]
		for all $z\in X$. Therefore, we obtain $\pi_{x,a,b}= \pi_{y,a,b}$, as claimed. Hence we may write $\pi_{a,b}=\pi_{x,a,b}$.
		Next, by Corollary~\ref{cor:basicFact}(1) we know that $\theta_{(x,a)}^2=\id$. This leads to
		\[(y,b)=\theta_{(x,a)}^2(y,b)=\theta_{(x,a)}(\pi_{a,b}(y),a+b)=(\pi_{a,a+b}\pi_{a,b}(y),a+(a+b))=(\pi_{a,a+b}\pi_{a,b}(y),b).\]
		Hence $\pi_{a,a+b}\pi_{a,b}=\id$ for all $a,b\in A$. Replacing $b$ by $a+b$ in the last equality we get
		$\pi_{a,b}\pi_{a,a+b}=\id$. Therefore, $\pi_{a,b}\in\Sym(X)$ and $\pi_{a,b}^{-1}=\pi_{a,a+b}$ for all $a,b\in A$.
		Furthermore, Corollary~\ref{cor:basicFact}(1) assures also that $\theta_{\theta_{(x,a)}(y,b)}=\theta_{(x,a)}\theta_{(y,b)}$,
		which yields
		\begin{align*}
			(\pi_{a+b,c}(z),(a+b)+c) & =\theta_{(\pi_{a,b}(y),a+b)}(z,c)=\theta_{\theta_{(x,a)}(y,b)}(z,c)\\
			& =\theta_{(x,a)}\theta_{(y,b)}(z,c)=\theta_{(x,a)}(\pi_{b,c}(z),b+c)\\
			& =(\pi_{a,b+c}\pi_{b,c}(z),a+(b+c)).
		\end{align*}
		Hence $\pi_{a,b+c}\pi_{b,c}=\pi_{a+b,c}$ for all $a,b,c\in A$. In particular, plugging $c=0$ to the last
		equality, we obtain $\pi_{a,b}=\pi_{a+b,0}\pi_{b,0}^{-1}$. Therefore, defining the map $\sigma\colon A\to\Sym(X)$
		by the formula $\sigma_a=\pi_{a,0}$, we conclude that $\theta_{(x,a)}(y,b)=(\sigma_{a+b}\sigma_b^{-1}(y),a+b)$
		and thus \[s((x,a),(y,b))=((x,a),(\sigma_{a+b}\sigma_b^{-1}(y),a+b)),\] which shows that $(S,s)=\Ext_X^\sigma(A,t_A)$.
	\end{proof}
\end{pro}

The following result shows that all solutions of the form $\Ext_X^\sigma(A,t_A)$ are, in fact, isomorphic.

\begin{pro}\label{pro:isom-of-sol-prime}
	Assume that $(A,+)$ is an elementary abelian $2$-groups and $X$ is a non-empty set. If $\sigma\colon A\to\Sym(X)$
	and $\rho\colon A\to\Sym(X)$ then the solutions $\Ext_X^\sigma(A,t_A)$ and $\Ext_X^\rho(A,t_A)$ of the PE are isomorphic.
	\begin{proof}
		Note first that if $\varphi\in\Sym(X)$ is an arbitrary permutation then defining $\pi\colon A\to\Sym(X)$ by the formula
		$\pi_a=\sigma_a\varphi$, we get $\pi_{a+b}\pi_b^{-1}=\sigma_{a+b}\sigma_b^{-1}$ for all $a,b\in A$. In particular, 
		$\Ext_X^\sigma(A,t_A)=\Ext_X^\pi(A,t_A)$. Hence we may assume that $\sigma_0=\id=\rho_0$. We claim that
		the map $f\colon\Ext_X^\sigma(A,t_A)\to\Ext_X^\rho(A,t_A)$, defined as $f(x,a)=(\rho_a\sigma_a^{-1}(x),a)$,
		is an isomorphism of solutions. Clearly, $f$ is a bijective map. Moreover, put $\Ext_X^\sigma(A,t_A)=(X\times A,s_\sigma)$
		and $\Ext_X^\rho(A,t_A)=(X\times A,s_\rho)$. Then
		\begin{align*}
			(f\times f)s_\sigma((x,a),(y,b)) & =(f\times f)((x,a),(\sigma_{a+b}\sigma_b^{-1}(y),a+b))\\
			& =((\rho_a\sigma_a^{-1}(x),a),(\rho_{a+b}\sigma_{a+b}^{-1}\sigma_{a+b}\sigma_b^{-1}(y),a+b))\\
			& =((\rho_a\sigma_a^{-1}(x),a),(\rho_{a+b}\sigma_b^{-1}(y),a+b))\\
			& =((\rho_a\sigma_a^{-1}(x),a),(\rho_{a+b}\rho_b^{-1}\rho_b\sigma_b^{-1}(y),a+b))\\
			& =s_\rho((\rho_a\sigma_a^{-1}(x),a),(\rho_b\sigma_b^{-1}(y),b))\\
			& =s_\rho(f\times f)((x,a),(y,b)).
		\end{align*}
		Therefore, $f$ is a morphism of solutions and the result is proved.
	\end{proof}
\end{pro}

Combining Theorem~\ref{thm:involuiveGxI} and Proposition~\ref{pro:retrsol-prime} we obtain a complete description
of all involutive solutions of the PE. Surprisingly, it turns out that all such solutions are uniquely determined,
up to isomorphism, by two elementary abelian $2$-groups.

\begin{thm}\label{thm:main-prime}
	Assume that $(S,s)$ is an involutive solution of the PE. Then there exist elementary abelian $2$-groups
	$A$ and $G$ and a non-empty set $X$ such that $S$ may be identified with $X\times A\times G$ and then
	\[(S,s)=\Ext_X^\sigma(A,t_A)\times(G,s_G)\] for some $\sigma\colon A\to\Sym(X)$, where $(A,t_A)$ is the
	unique irretractable involutive solution of the PE on $A$ and $(G,s_G)$ is the unique bijective solution of the PE on $G$.
	Moreover, $\Ret(S,s)=(A,t_A)$ and $(S,s)$ is isomorphic to the product $(X,\id)\times(A,t_A)\times(G,s_G)$ of solutions.
	\begin{proof}
		In view of what was already proved, it remains to show that the solution $(S,s)$ is isomorphic to
		the product $(X,\id)\times(A,t_A)\times(G,s_G)$. But Proposition~\ref{pro:isom-of-sol-prime}
		yields $\Ext_X^\sigma(A,t_A)\cong\Ext_X^\rho(A,t_A)$, where $\rho\colon A\to\Sym(X)$ is defined as $\rho_a=\id$
		for all $a\in A$. Since we clearly have $\Ext_X^\rho(A,t_A)\cong(X,\id)\times(A,t_A)$, the result is proved.
	\end{proof}
\end{thm}

\begin{thm}\label{thm:main2}
	All involutive solutions, up to isomorphism, of the PE defined on a non-empty set $S$ are in a bijective correspondence
	with decompositions of $S$ as a product $X\times A\times G$, where $X$ is a non-empty set and $A,G$ are elementary abelian $2$-groups.
	\begin{proof}
		Let $(S_1,s_1)$ and $(S_2,s_2)$ be involutive solutions of the PE corresponding to decompositions
		of $S$ of the form $S_1=X_1\times A_1\times G_1$ and $S_2=X_2\times A_2\times G_2$, respectively (see Theorem~\ref{thm:main-prime}).
		
		If $|X_1|=|X_2|$, $|A_1|=|A_2|$ and $|G_1|=|G_2|$ then there exist a bijection $\alpha\colon X_1\to X_2$ and morphisms of groups
		$\beta\colon(A_1,+)\to(A_2,+)$ and $\gamma\colon(G_1,\cdotp)\to(G_2,\cdotp)$ (see the proof of Proposition~\ref{dim-arg}).
		But then it is easy to check that the map $f\colon(S_1,s_1)\to(S_2,s_2)$, defined as $f(x,a,g)=(\alpha(x),\beta(a),\gamma(g))$
		for $x\in X_1$, $a\in A_1$ and $g\in G_1$, is an isomorphism of solutions.
		
		Conversely, suppose that $f\colon(S_1,s_1)\to(S_2,s_2)$ is an isomorphism of solutions. Then \eqref{eq:hom} assures that
		$f\colon S_1\to S_2$ is a morphism of semigroups (recall that multiplication in $S_i$ for $1\leq i\leq 2$ is given as
		$(x,a,g)\cdot(y,b,h)=(x,a,g\cdot h)$ for $x,y\in X_i$, $a,b\in A_i$ and $g,h\in G_i$). Since $G_i=e_iS_ie_i$ for
		any $e_i\in E(S_i)$ (see the comment before Corollary~\ref{cor:basicFact}), we get $f(G_1)=f(eS_1e)=f(e)S_2f(e)=G_2$
		for each $e\in E(S_1)$. Hence $|G_1|=|G_2|$. Next, by Theorem~\ref{thm:main-prime}, we obtain
		\[(A_1,t_{A_1})=\Ret(S_1,s_1)\cong\Ret(S_2,s_2)=(A_2,t_{A_2})\] and thus $|A_1|=|A_2|$. Finally, to prove that $|X_1|=|X_2|$,
		write $f(x,0,1)=(\varphi(x),\psi(x),1)$ for $x\in X_1$, where $\varphi\colon X_1\to X_2$ and $\psi\colon X_1\to A_2$
		(we know that $f(x,0,1)$ is of that form because $(x,0,1)\in E(S_1)$ implies $f(x,0,1)\in E(S_2)=X_2\times A_2\times\{1\}$). Then
		\begin{align*}
			((\varphi(x),\psi(x),1),(\varphi(y),\psi(y),1)) & =(f\times f)s_1((x,0,1),(y,0,1))\\
			& =s_2(f\times f)((x,0,1),(y,0,1))\\
			& =s_2((\varphi(x),\psi(x),1),(\varphi(y),\psi(y),1)\\
			& =((\varphi(x),\psi(x),1),(\varphi(y),\psi(x)+\psi(y),1)
		\end{align*}
		for all $x,y\in X_1$. This implies that $\psi(x)=0$ for all $x\in X_1$.
		Hence $f(x,0,1)=(\varphi(x),0,1)$. Since the inverse $f^{-1}$ of $f$ has a similar property,
		we obtain that $\varphi$ is a bijection and thus $|X_1|=|X_2|$. Hence the proof is finished.
	\end{proof}
\end{thm}

As an immediate consequence of Theorem~\ref{thm:main2} we get the following result.

\begin{cor}\label{cor:binom}
	Assume that $S$ is a finite set of cardinality $|S|=2^n(2m+1)$ for some $n,m\geq 0$. Then there exist,
	up to isomorphism, exactly $\binom{n+2}{2}$ involutive solutions of the PE defined on $S$.
\end{cor}

As an illustration of the main results of this section we provide the form of all six, up to isomorphism,
involutive solutions of the PE defined on a set of cardinality $12$.

\begin{exa}
	Let $S$ be a set of cardinality $12$. The involutive solutions $(S,s)$ of the PE presented below correspond to the decomposition
	$S=X\times A\times G$, where $(A,+)$ and $(G,\cdotp)$ are elementary abelian $2$-groups. Clearly, we have
	$A,G\in\{0,C_2,C_2\times C_2\}$. Moreover, if one of the factors $A$ or $G$ is the trivial group $0$ then this factor
	will be omitted and thus we write $S=X$ if $A=G=0$ or $S=X\times A$ if $A\ne 0$ but $G=0$ or $S=X\times G$ if $A=0$ but $G\ne 0$.
	\begin{enumerate}
		\item Let $X=\{1,\dotsc,12\}$. Then $s(x,y)=(x,y)$.
		\item Let $X=\{1,\dotsc,6\}$ and $A=C_2$. Then $s((x,a),(y,b))=((x,a),(y,a+b))$.
		\item Let $X=\{1,\dotsc,6\}$ and $G=C_2$. Then $s((x,g),(y,h))=((x,g\cdot h),(y,h))$.
		\item Let $X=\{1,2,3\}$ and $A=C_2\times C_2$. Then $s((x,a),(y,b))=((x,a),(y,a+b))$.
		\item Let $X=\{1,2,3\}$ and $G=C_2\times C_2$. Then $s((x,g),(y,h))=((x,g\cdot h),(y,h))$.
		\item Let $X=\{1,2,3\}$ and $A=G=C_2$. Then $s((x,a,g),(y,b,h))=((x,a,g\cdot h),(y,a+b,h))$.
	\end{enumerate}
\end{exa}

\section{Structure monoids and algebras of solutions of the PE}

Suppose $S$ is a non-empty finite set and $s\colon S\times S\to S\times S$ is a bijective map. One simply says that $(S,s)$
is a quadratic set. Then, over a field $K$, one associates a quadratic algebra, called the structure algebra of $(S,s)$ over $K$,
\[A(K,S,s)=K\free{S \mid xy=uv\text{ if }s(x,y)=(u,v)}.\] This algebra has a natural gradation $A(K,S,s)=\bigoplus_{n\geq 0}A_n$
and is generated by $A_1=\Span_K(S)$. So $A(K,S,s)$ is a connected graded $K$-algebra and each $A_n$ is finite-dimensional. Quadratic
algebras have received a lot of attention and for more information we refer to \cite{PolPos}. Again because the defining relations
are homogeneous, \[A(K,S,s)=K[M(S,s)],\] the monoid $K$-algebra of the so-called structure monoid $M(S,s)$ of $(S,s)$ defined as
\[M(S,s)=\free{S \mid xy=uv\text{ if }s(x,y)=(u,v)}.\]

In the last two decades there has been an intense study of monoids, groups and algebras associated with quadratic sets,
see for example \cite{CedOkn,GI2012,GI2018,JesOknBook}. In case $(S,s)$ is a finite non-degenerate unitary solution of the
Quantum Yang--Baxter Equation, i.e., $(S,\tau s)$ is a finite non-degenerate involutive solution of the Yang--Baxter Equation
then it has been shown in \cite{GIVdB} that the structure algebra $A(K,S,\tau s)$ has a very rich structure that shares many
properties with polynomial algebras in finitely many commuting variables. In particular, it is a (left and right) Noetherian
algebra that satisfies a polynomial identity (abbreviated, PI) and, moreover, it is a domain that has finite Gelfand--Kirillov
dimension. Hence, $M(S,\tau s)$ has a group of quotients that is torsion-free and it is abelian-by-finite. The latter also has
been investigated in \cite{ESS}. Note that if $s=\id$ then $K[M(S,s)]=K\free{X}$, the free $K$-algebra in variables
in the set $X=\{x_s:x\in S\}$, while $K[M(S,\tau s)]=K[X]$, the polynomial algebra in commuting variables in $X$.
In case $(S,s)$ is an arbitrary finite, bijective and non-degenerate solution of the Yang--Baxter Equation, it has been
shown in \cite{JesKubVan,JesKubVan2} that the algebraic structure of $A(K,S,\tau s)$ determines when $s$ is involutive.

In this section we investigate the structure algebra $A(K,S,\tau s)=K[M(S,\tau s)]$ corresponding to an involutive finite solution $(S,s)$
of the PE. The description of the solutions obtained in the previous section allows to show that $M(S,\tau s)$ is a finite extension of a
free abelian submonoid and hence one obtains that $K[M(S,\tau s)]$ is a Noetherian PI-algebra.

Let $(S,s)$ be a solution of the PE. So, as before, we write $s(s,y)= (xy,\theta_{x}(y))$. To avoid possible confusion between
the multiplication in the semigroup $S$ and the product in $M(S,\tau s)$, we write the product in $M(S,\tau s)$ as $\circ$.
Thus \[M(S,\tau s)=\free{S\mid x\circ y=\theta_x(y)\circ xy\text{ for all }x,y\in S}.\] For $x\in S$ and $n\geq 1$
we shall write $x^{(n)}=x\circ\dotsm\circ x$ (the product of $n$ copies of $x$ in $M(S,\tau s)$) to distinguish this element
from the power $x^n$ of $x$ in $S$. Moreover, we put $x^{(0)}=1$, the identity element of $M(S,\tau s)$.

In order to state the main result of this section we recall some well-known notation. The rank $\rk S$ of a semigroup $S$
is defined as the supremum of the ranks of free abelian subsemigroups of $S$. The classical Krull dimension of an algebra $A$
is denoted by $\clK A$. Moreover, $\GK A$ stands for the Gelfand--Kirillov dimension of $A$. For more information and background
on ring theory we refer the reader to \cite{MR}.

\begin{thm}\label{algrep}
	Assume that $(S,s)$ is a finite involutive solution of the PE. Then the structure monoid $M=M(S,\tau s)$ is (free abelian)-by-finite.
	More precisely, there exist a free abelian submonoid $C\subseteq M$ of rank $n=|E(S)|\cdot|{\Ret(S,s)}|^{-1}$ and finite subsets
	$L,R\subseteq M$ such that $M=\bigcup_{f\in L}C\circ f=\bigcup_{f\in R}f\circ C$. In particular, if $K$ is a field then the structure
	algebra $A(K,S,\tau s)=K[M]$ is a module finite extension of the polynomial algebra $K[C]$ in $n$ commuting variables. Hence $K[M]$ is
	a Noetherian PI-algebra satisfying \[\clK K[M]=\GK K[M]=\rk M=n\leq|S|\] and the equality holds if and only if $s=\id$.
	\begin{proof}
		By Theorem~\ref{thm:main-prime} we may assume that $S=X\times A\times G$, where $X$ is a non-empty set,
		$(A,+)$ and $(G,\cdotp)$ are elementary abelian $2$-groups and \[s((x,a,g),(y,b,h))=((x,a,g\cdot h),(y,a+b,h)).\]
		Then $S$ is a semigroup with multiplication given by the formula $(x,a,g)\cdot(y,b,h)=(x,a,g\cdot h)$. Since
		$E(S)=X\times A\times\{1\}$ and $\Ret(S,s)=(A,t_A)$ (see Theorem~\ref{thm:main-prime}), we get $n=|E(S)|\cdot|{\Ret(S,s)}|^{-1}=|X|$.
		
		Write $X=\{x_1,\dotsc,x_n\}$. Since \[(x,a,g)\circ(y,b,h)=(y,a+b,h)\circ(x,a,g\cdot h)\] for all $x,y\in X$, $a,b\in A$
		and $g,h\in G$, it follows that each element $w\in M$ may be written in the following form $w=w_1\circ\dotsb\circ w_n$, where
		$w_i\in\free{(x_i,a,g):a\in A\text{ and }g\in G}\subseteq M$ for $1\leq i\leq n$. Moreover, writing $A=\{a_1,\dotsc,a_m\}$
		(clearly $m=2^r$ for some $r\geq 0$) and using the fact that \[(x,a,g)\circ(x,b,h)=(x,a+b,h)\circ(x,a,g\cdot h)=(x,b,h)\circ(x,a+b,g)\]
		for all $x\in A$, $a,b\in A$ and $g,h\in G$, we obtain that each $w_i$ may be written as $w_i=w_{i1}\circ\dotsb\circ w_{im}$,
		where $w_{ij}\in\free{(x_i,a_j,g):g\in G}$ for $1\leq j\leq m$. Finally, the identity
		\[(x,a,g)\circ(x,a,h)=(x,0,h)\circ(x,a,g\cdot h),\] holding for all $x\in X$, $a\in A$ and $g,h\in G$, yields
		$w_{ij}\in\free{(x_i,a_j,g):g\in G}\subseteq\bigcup_{g\in G}N\circ(x_i,a_j,g)$, where
		\[N=\free{(x,0,g):x\in X\text{ and }g\in G}\subseteq M.\]
		Furthermore, we have \[(x,0,g)\circ(y,b,h)=(y,b,h)\circ(x,0,g\cdot h)\] for all $x,y\in X$, $b\in A$ and $g,h\in G$.
		Replacing $g$ by $g\cdot h$ in the above formula, we get \[(y,b,h)\circ(x,0,g)=(x,0,g\cdot h)\circ(y,b,h)\] as well.
		Hence $w\circ N=N\circ w$ for each $w\in M$, that is, $N$ is a normal submonoid of $M$. Therefore,
		\begin{align*}
			w_i & =w_{i1}\circ\dotsb\circ w_{im}\in\bigcup_{\substack{g_1,\dotsc,g_m\in G\\ \varepsilon_1,\dotsc,\varepsilon_m\in\{0,1\}}}
			N\circ(x_i,a_1,g_1)^{(\varepsilon_1)}\circ\dotsb\circ(x_i,a_m,g_m)^{(\varepsilon_m)},
		\end{align*}
		which leads to $w=w_1\circ\dotsb\circ w_n\in\bigcup_{f\in F_1}N\circ f$, where $F_1\subseteq M$ is a finite subset of $M$
		consisting of elements of the form
		\[(x_1,a_1,g_{11})^{(\varepsilon_{11})}\circ\dotsb\circ(x_1,a_m,g_{1m})^{(\varepsilon_{1m})}
		\circ\dotsb\circ(x_n,a_1,g_{n1})^{(\varepsilon_{n1})}\circ\dotsb\circ(x_n,a_m,g_{nm})^{(\varepsilon_{nm})},\]
		where $g_{ij}\in G$ and $\varepsilon_{ij}\in\{0,1\}$ for $1\leq i\leq n$ and $1\leq j\leq m$
		(in particular, $|F_1|\leq(2|G|)^{nm}$). 
		Thus we have proved that
		\begin{equation}\label{sum-1}
		M=\bigcup_{f\in F_1}N\circ f=\bigcup_{f\in F_1}f\circ N.
		\end{equation}
		Next, write $G=\{g_1,\dotsc,g_k\}$ (clearly $k=2^t$ for some $t\geq 0$). Since $w\circ N=N\circ w$ for each
		$w\in M$, we get that $N$ consists of elements of the form
		\[(x_1,0,g_1)^{(l_{11})}\circ\dotsb\circ(x_1,0,g_k)^{(l_{1k})}\circ\dotsb\circ(x_n,0,g_1)^{(l_{n1})}\circ\dotsb\circ(x_n,0,g_k)^{(l_{nk})},\]
		where $l_{ij}\geq 0$ for $1\leq i\leq n$ and $1\leq j\leq k$. But it is easy to see that $(x,0,g)^{(l+1)}=(x,0,1)^{(l)}\circ(x,0,g)$
		for each $l\geq 0$. Since $(x,0,g)\circ(y,0,1)=(y,0,1)\circ(x,0,g)$ for all $x,y\in X$ and $g\in G$, it follows that
		$N=\bigcup_{f\in F_2}C\circ f$, where \[C=\free{(x,0,1):x\in X}\subseteq\Cen(N),\]
		$\Cen(N)$ is the center of $N$, and $F_2\subseteq N$ is a finite subset consisting of elements of the form
		\[(x_1,0,g_1)^{(\varepsilon_{11})}\circ\dotsb\circ(x_1,0,g_k)^{(\varepsilon_{1k})}\circ\dotsb\circ
		(x_n,0,g_1)^{(\varepsilon_{n1})}\circ\dotsb\circ(x_n,0,g_k)^{(\varepsilon_{nk})},\]
		where $\varepsilon_{ij}\in\{0,1\}$ for $1\leq i\leq n$ and $1\leq j\leq k$ (in particular $|F_2|\leq 2^{nk}$). Therefore, we have
		\begin{equation}\label{sum-2}
		N=\bigcup_{f\in F_2}C\circ f=\bigcup_{f\in F_2}f\circ C.
		\end{equation}
		Combining \eqref{sum-1} and \eqref{sum-2} we conclude that $M=\bigcup_{f\in L}C\circ f=\bigcup_{f\in R}f\circ C$, where
		\[L=\{f_2\circ f_1:f_1\in F_1\text{ and }f_2\in F_2\},\qquad R=\{f_1\circ f_2:f_1\in F_1\text{ and }f_2\in F_2\}.\]
		Clearly, $L,R\subseteq M$ are finite subsets of $M$. Moreover, because $(x,0,1)\circ(y,0,1)=(y,0,1)\circ(x,0,1)$
		for all $x,y\in X$, we conclude that $C$ is a free abelian monoid of rank $|X|=n$, which finishes the proof of first
		part of the theorem.
		
		Since $K[M]=\sum_{f\in L}K[C]\circ f=\sum_{f\in R}f\circ K[C]$ is a module finite extension of $K[C]\cong K[t_1,\dotsc,t_n]$,
		the polynomial algebra in $n$ commuting variables, it follows that $K[M]$ is a Noetherian PI-algebra (see \cite{MR}).
		Hence, a result of Ananin \cite{Ana1989} implies that $M$ is a linear monoid and then
		\cite[Proposition~1, p.~221, Proposition~7, p.~280--281, and Theorem~14, p.~284]{Okn1991} yield
		\[\clK K[M]=\GK K[M]=\rk M=n\leq|S|.\]
		
		Clearly, we have equality in the above if and only if $S=X$ (i.e., $A=G=0$, the trivial group). Hence, the last
		part of the theorem follows as well.
	\end{proof}
\end{thm}

Note that if $(S,s)$ is an infinite solution of the PE then the structure algebra $A=A(K,S,\tau s)$ is not
Noetherian as the ideal $I=\sum_{x\in S}AxA$ is not finitely generated neither as a left nor right ideal.

As an illustration  of Theorem~\ref{algrep} we give two applications in case $(S,s)$ is a finite involutive solution
of the PE based on a left zero semigroup (i.e., $G=0$; see Theorem~\ref{thm:main2}). First we show that the prime ideals
of the structure algebra $A=A(K,S,\tau s)$ over a field $K$ are very special and second we show that certain ring-theoretical
properties of $A$ encode information about the map $s$.

\begin{pro}\label{alg-1}
	Assume that $(S,s)$ is a finite involutive solution of the PE, where $S$ is a left zero semigroup. Let $M=M(S,\tau s)$.
	If $P$ is a prime ideal of the structure algebra $A(K,S,\tau s)=K[M]$ over a field $K$ then the algebra $K[M]/P$ is a commutative
	affine domain. In particular, each prime ideal of $K[M]$ is completely prime and semiprimitive and $K[M]$ is a Jacobson ring.
	\begin{proof}
		We use the same notation as in Theorem~\ref{algrep} and its proof. Because $S$ is a left zero semigroup, it easily
		is verified that each element $w$ of $M$ is normal in $M$, i.e., $w\circ M=M\circ w$. Let $P$ be a prime ideal of $K[M]$.
		Since $x\circ x=\theta_x(x)\circ x$ for each $x\in S$, we get $(x-\theta_x(x))\circ x=0$ in $K[M]$. Thus, by normality of
		elements of $M$, it follows that $(x-\theta_x(x))\circ K[M]\circ x=0$. This yields $x-\theta_x(x)\in P$ or $x\in P$.
		Define $X=S\cap P$ and $Y=S\setminus X$. By what we have already shown, it follows that $y-\theta_y(y)\in P$ for each $y\in Y$.
		Moreover, if $x,y\in Y$ then $\theta_x(y)\in Y$. Indeed, otherwise $x\circ y=\theta_x(y)\circ x\in P$ and thus, again by normality
		of elements of $M$, we get $x\in P$ or $y\in P$, a contradiction. Since $\theta_y(y)\circ x=x\circ\theta_y(y)$, we obtain 
		$(\theta_x(y)-\theta_y(y))\circ x=x\circ(y-\theta_y(y))\in P$, which implies $\theta_x(y)-\theta_y(y)\in P$ because $x\notin P$.
		Thus, in consequence, $y-\theta_x(y)=(y-\theta_y(y))-(\theta_x(y)-\theta_y(y))\in P$ for all $x,y\in Y$.
		Therefore, if $P_0$ is the ideal of $K[M]$ generated by $X$ and elements of the form $y-\theta_x(y)$ for all $x,y\in Y$
		then $P_0\subseteq P$ and
		\[K[M]/P_0\cong\frac{K\free{Y\mid x\circ y=\theta_x(y)\circ x\text{ for all }x,y\in Y}}{(y-\theta_x(y):x,y\in Y)}
		\cong K[t_1,\dotsc,t_k],\]
		the polynomial algebra in $k$ commuting variables, where $k$ is equal to the number of orbits of the action of the
		group $\{\theta_x:x\in Y\}\subseteq\Sym(Y)$ on the set $Y$. In particular, the ideal $P_0$ is prime and semiprimitive
		(note also that if $P$ is a minimal prime ideal then $P=P_0$ and we have a description of $P$ in terms of generators).
		Since $K[M]/P$ is a homomorphic image of the commutative affine algebra $K[M]/P_0$, it follows that the ideal $P$
		is completely prime and semiprimitive (by Nullstellensatz, see \cite{MR}). Hence the fact that $K[M]$ is a Jacobson ring
		also follows.
	\end{proof}
\end{pro}

It is worth to add that (under the assumptions and notation from Proposition~\ref{alg-1}) the Jacobson radical $J=\mathcal{J}(K[M])$ of $K[M]$
coincides with the prime radical $\mathcal{B}(K[M])$ of $K[M]$ and that it is a nilpotent ideal, say of index $n$. But then the fact that the
algebra $K[M]/J$ is commutative leads to a conclusion that $(xy-yx)^n$ is a polynomial identity satisfied by $A(K,S,\tau s)=K[M]$.

As a consequence of Proposition~\ref{alg-1} we get the following result.

\begin{cor}\label{alg-2}
	Assume that $(S,s)$ is a finite involutive solution of the PE, where $S$ is a left zero semigroup. Let $M=M(S,\tau s)$.
	If $K$ is a field then the structure algebra $A(K,S,\tau s)=K[M]$ is a domain if and only if it is prime if and only if $s=\id$.
	\begin{proof}
		Clearly $s=\id$ implies that $K[M]$ is a domain (because in this case $K[M]\cong K[t_1,\dotsc,t_n]$, the polynomial algebra
		in $n=|S|$ commuting variables), and if $K[M]$ is a domain then $K[M]$ is prime. So, it remains to show that
		primeness of $K[M]$ implies that $s=\id$ or, equivalently, that $\theta_x=\id$ for each $x\in S$. But if $K[M]$ is prime,
		Proposition~\ref{alg-1} implies that $K[M]$ is a commutative domain. In particular, $M$ is a commutative and cancellative monoid.
		But then $x\circ y=\theta_x(y)\circ x=x\circ\theta_x(y)$ for $x,y\in S$ implies $\theta_x(y)=y$ in $M$. Since the defining
		relations of $M$ are quadratic, it follows that $S$ embeds in $M$ and thus $\theta_x(y)=y$ in $S$, which finishes the proof.
	\end{proof}
\end{cor}

We finish with an example that shows that the previous corollary cannot be extended to the semiprime case.
Indeed, we give an example that is semiprime (or equivalently semiprimitive) but not a domain.

\begin{exa}
	Let $S=C_2$ be the additively written cyclic group of order two and $s(x,y)=(x,x+y)$ for $x,y\in S$. Then
	\[M=M(S,\tau s)=\free{x,y\mid x\circ x=x\circ y=y\circ x}\] and thus, for a field $K$, we get $A(K,S,\tau s)=K[M]\cong K[u,v]/(u^2-uv)$.
	Denoting the images of variables $u,v$ in $R=K[u,v]/(u^2-uv)$ again by $u,v$, we get $R=K[v]\oplus uK[v]$. Suppose that $r\in R$
	satisfies $r^2=0$. Write $r=a+ub\in R$ for some $a,b\in K[v]$. Then $0=r^2=a^2+u(2ab+b^2v)$ implies $a^2=0$ and $2ab+b^2v=0$.
	Thus $a=0$ and then $b^2v=0$, which leads to $b=0$ (because $K[v]$, as a polynomial algebra, is a domain). Therefore,
	$r=0$ and thus the algebra $R\cong K[M]$ is reduced. Since $K[M]$ is commutative, it is semiprime but not a domain.
\end{exa}

\end{document}